\documentclass{article}

\usepackage{amssymb}
\usepackage{amsmath}

\usepackage[dvips]{graphicx}

\begin{document}


\begin{center}
\textbf{ A note on $n$-Jordan homomorphisms }
\end{center}

\bigskip

\begin{center}
M. El Azhari 
\end{center}

\bigskip

\noindent\textbf{Abstract.} By using a variation of a theorem on $n$-Jordan homomorphisms due to Herstein, we deduce the following G. An's result: Let $ A $ and $ B $ be two rings where $ A $ has a unit and $ char(B)> n. $ If every Jordan homomorphism from $ A $ into $ B $ is a homomorphism  (anti-homomorphism), then every $n$-Jordan homomorphism from $ A $ into $ B $ is an $n$-homomorphism (anti-$n$-homomorphism).
\\
\\
\noindent \textbf{1. Preliminaries }\\
\\
\noindent\textbf{} Let $ A, B $  be two rings and $ n\geqslant 2 $ an integer. An additive map $ h:A\rightarrow B $ is called an $n$-Jordan homomorphism if $ h(x^{n})=h(x)^{n} $ for all  $ x\in A. $ Also, an additive map $ h:A\rightarrow B $  is called an $n$-homomorphism or an anti-$n$-homomorphism if $ h(\prod_{i=1}^{n}x_{i})=\prod_{i=1}^{n} h(x_{i}) $ or $ h(\prod_{i=1}^{n}x_{i})=\prod_{i=0}^{n-1} h(x_{n-i}), $ respectively, for all  $ x_{1},...,x_{n}\in A. $ \\
In the usual sense, a $2$-Jordan homomorphism is a Jordan homomorphism, a $2$-homomorphism is a homomorphism and an anti-$2$-homomorphism is an anti-homomorphism.
It is obvious that $n$-homomorphisms are $n$-Jordan homomorphisms. Conversely, under certain conditions, $n$-Jordan homomorphisms are 
$n$-homomorphisms.\\
We say that a ring $ A $ is of characteristic greater than $ n\:\:(char(B)>n) $ if $ n!x=0 $ implies $ x=0 $ for all $ x\in A. $\\
\\
\noindent \textbf{2. Results}\\
\\
\noindent \textbf{ Lemma 2.1} [4, Lemma 1]. Let  $ A,\:B $ be two rings, $ n\geqslant 2 $ be an integer and  $ f:A^{n}\rightarrow B $ be a multi-additive map such that $ f(x,x,...,x)=0 $ for all $ x $ in $ A. $  Then $ \: \sum_{\sigma\in S_{n}} f(x_{\sigma(1)},...,x_{\sigma(n)})=0 $ for all $ x_{1},...,x_{n}\in A, $  where $ S_{n} $ is the set of all permutations on $ \{1,...,n\}.$\\
\\
\noindent \textbf{} By using Lemma 2.1, we have the following lemma:\\
\\
\noindent \textbf{ Lemma 2.2.} Let  $ A,\:B $ be two rings, $ n\geqslant 2 $ be an integer and $ h:A\rightarrow B $ be an $n$-Jordan homomorphism.  Then \\
$ \sum_{\sigma\in S_{n}} (h(\prod_{i=1}^{n}x_{\sigma(i)})-\prod_{i=1}^{n} h(x_{\sigma(i)}))=0 $ for all $ x_{1},...,x_{n}\in A, $
where $ S_{n} $ is the set of all permutations on $ \{1,...,n\}. $ \\
\\
\noindent \textbf{Proof.} Consider the map $ f:A^{n}\rightarrow B,$  \\
$ f(x_{1},x_{2},...,x_{n})= h(\prod_{i=1}^{n}x_{i})-\prod_{i=1}^{n} h(x_{i}),\: f $ is clearly multi-additive and\\
$  f(x,x,...,x)= h(x^{n})-h(x)^{n}=0 $ for all $ x\in A. $ By Lemma 2.1,\\
$  \sum_{\sigma\in S_{n}} f(x_{\sigma(1)},...,x_{\sigma(n)})= \sum_{\sigma\in S_{n}} (h(\prod_{i=1}^{n}x_{\sigma(i)})-\prod_{i=1}^{n} h(x_{\sigma(i)}))= 0 $ for all  $ x_{1},...,x_{n}\in A. $\\
\\ 
It was shown in [3] that if $ n\geqslant 2,\: A,\: B $ are commutative rings, $ char(B)>n $ and  $ h:A\rightarrow B $ is an $n$-Jordan homomorphism, then $ h $ is an $n$-homomorphism. The same was also proved for algebras in [2] and [5]. Here we obtain this result as a  consequence of Lemma 2.2.\\
\\
\noindent \textbf{ Theorem 2.3.} Let  $ A,\:B $ be two commutative rings and $ char(B)>n. $ Then every $n$-Jordan homomorphism from $ A $ into $ B $ is an $n$-homomorphism.\\
\\
\noindent \textbf{Proof.} Let $ h:A\rightarrow B $ be an $n$-Jordan homomorphism. By Lemma 2.2 and since $ A,\: B $ are commutative,\\
$ \sum_{\sigma\in S_{n}} (h(\prod_{i=1}^{n}x_{\sigma(i)})-\prod_{i=1}^{n} h(x_{\sigma(i)}))= n! (h(\prod_{i=1}^{n}x_{i})-\prod_{i=1}^{n} h(x_{i}))= 0 $ for all $ x_{1},...,x_{n}\in A, $ hence $ h $ is an $n$-homomorphism since $ char(B)>n. $\\
\\
Now we give the following variation of a theorem on $n$-Jordan homomorphisms due to Herstein [4, Theorem K].\\
\\
\noindent \textbf{ Theorem 2.4.} Let $ h $ be an $n$-Jordan homomorphism from a ring $ A $ into a ring $ B $ with $ char(B)>n. $ Suppose further that $ A $ has a unit $ e, $ then $ h=h(e)\tau $ where $ h(e) $ is in the centralizer of $ h(A) $ and $ \tau $ is a Jordan homomorphism.\\
\\
\noindent \textbf{Proof.} Since $ h $ is an $n$-Jordan homomorphism, $ h(e)=h(e^{n})=h(e)^{n}. $ By Lemma 2.2 and putting $ x_{1}=x,\:x_{2}=x_{3}=\cdots=x_{n}=e, $  \\ 
$ n!h(x)=(n-1)!(h(e)^{n-1}h(x)+h(e)^{n-2}h(x)h(e)+\dots+h(x)h(e)^{n-1}) $ \\
and so $ nh(x)=h(e)^{n-1}h(x)+h(e)^{n-2}h(x)h(e)+\dots+h(x)h(e)^{n-1}\:\:\: (1) $\\
since $ char(B)>n. $\\
By multiplying on the right by $ f(e) $ both sides of the equality (1),\\
$ nh(x)h(e)=h(e)^{n-1}h(x)h(e)+h(e)^{n-2}h(x)h(e)^{2}+\dots+h(x)h(e)\:\:\: (2) $\\
Also, by multiplying on the left by $ f(e) $ both sides of the equality (1),\\
$ nh(e)h(x)=h(e)h(x)+h(e)^{n-1}h(x)h(e)+\dots+h(e)h(x)h(e)^{n-1}\:\:\: (3) $\\
By (2) and (3), $ (n-1)h(x)h(e)=(n-1)h(e)h(x)$ and consequently \\
$ h(x)h(e)=h(e)h(x)\:\:\: (4) $ since $ char(B)>n. $ Then $ h(e) $ is in the centralizer of $ h(A) .$
By (1) and (4), $ nh(x)=nh(e)^{n-1}h(x) $ and so $ h(x)=h(e)^{n-1}h(x)\:\:\: (5) $ since $ char(B)>n. $
By Lemma 2.2, (4) and putting  $ x_{1}=x_{2}=x,\:x_{3}=\cdots=x_{n}=e, \: n!(h(x^{2})-h(e)^{n-2}h(x)^{2})=0 $ and hence   $h(x^{2})=h(e)^{n-2}h(x)^{2}\:\:\: (6) $  since $ char(B)>n. $ 
Consider the map  $ \tau:A\rightarrow B,\:  \tau(x)=h(e)^{n-2}h(x),\: \tau $ is clearly additive. By (5), $ h(x)=h(e)^{n-1}h(x)=h(e)h(e)^{n-2}h(x)=h(e)\tau(x) $ for all $ x\in A. $ By (6),
$  \tau(x^{2})=h(e)^{n-2}h(x^{2})=h(e)^{2(n-2)}h(x)^{2}=(h(e)^{(n-2)}h(x))^{2}=\tau(x)^{2} $ for all $ x\in A, $ thus $ \tau $ is a Jordan homomorphism.\\
\\
As a consequence, we obtain the following result of G. An [1].\\
\\
\noindent \textbf{ Corollary 2.5} [1, Theorem 2.4]. Let $ A $ and $ B $ be two rings where $ A $ has a unit $ e $ and $ char(B)> n. $ If every Jordan homomorphism from $ A $ into $ B $ is a homomorphism  (anti-homomorphism), then every $n$-Jordan homomorphism from $ A $ into $ B $ is an $n$-homomorphism (anti-$n$-homomorphism).\\
\\
\noindent \textbf{Proof.} Let $ h:A\rightarrow B $ be an $n$-Jordan homomorphism. By Theorem 2.4, $ h=h(e)\tau $ where $ h(e) $ is in the centralizer of $ h(A) $ and $ \tau $ is a Jordan homomorphism. If $ \tau $ is a homomorphism, \\
$ h(x_{1}x_{2}\cdots x_{n})= h(e)\tau(x_{1}x_{2}\cdots x_{n})\\
=h(e)\tau(x_{1})\tau(x_{2})\cdots \tau(x_{n})\\
= h(e)^{n}\tau(x_{1})\tau(x_{2})\cdots \tau(x_{n})\:$ since $ h(e)=h(e^{n})=h(e)^{n} \\
= h(e)\tau(x_{1})h(e)\tau(x_{2})\cdots h(e)\tau(x_{n})\: $ since $ h(e) $ commutes with each $ \tau(x)\\
=h(x_{1})h(x_{2})\cdots h(x_{n}) $ for all $ x_{1},...,x_{n}\in A, $\\
hence $ h $ is an $n$-homomorphism. Similarly, if $ \tau $ is an anti-homomorphism, then $ h $ is an anti-$n$-homomorphism.

\bigskip

\begin{center}
REFERENCES
\end{center}

\bigskip

\noindent \textbf{} [1] G. An, Characterizations of $n$-Jordan homomorphisms, Linear and Multilinear Algebra, 66(4)(2018), 671-680.
 
\noindent \textbf{} [2] A. Bodaghi and H. Inceboz, $n$-Jordan homomorphisms on commutative algebras, Acta Math. Univ. Comenianae, 87(1)(2018), 141-146.

\noindent \textbf{} [3] E. Gselmann, On approximate $n$-Jordan homomorphisms, Annales Mathematicae Silesianae, 28(2014), 47-58.

\noindent \textbf{} [4] I. N. Herstein, Jordan homomorphisms, Trans. Amer. Math. Soc., 81(2)(1956), 331-341.

\noindent \textbf{} [5] Yang-Hi Lee, Stability of $n$-Jordan homomorphisms from a normed algebra to a Banach algebra, Abstract and Applied Analysis, 2013(2013), Article ID 691025, 5 pages.

\bigskip

\noindent \textbf{} Ecole Normale Sup\'{e}rieure

\noindent \textbf{} Avenue Oued Akreuch

\noindent \textbf{} Takaddoum, BP 5118, Rabat

\noindent \textbf{} Morocco
 
\bigskip 

\noindent \textbf{} E-mail:  mohammed.elazhari@yahoo.fr

\end{document}